\documentclass[submission,copyright,creativecommons]{eptcs}
\providecommand{\event}{SOS 2007} 

\usepackage{iftex}
\usepackage{amsthm,amsmath,amssymb}

\ifpdf
  \usepackage{underscore}         
  \usepackage[T1]{fontenc}        
\else
  \usepackage{breakurl}           
\fi

\newtheorem{lem}{Lemma}[section]
\newtheorem{thm}[lem]{Theorem}

\numberwithin{equation}{section}

\newcommand{\DDD}{{\mathcal{D}}}

\newcommand{\inv}{{\textsf{inv}}}

\newcommand{\exc}{{\textsf{exc}}}

\newcommand{\depth}{{\textsf{depth}}}

\newcommand{\fix}{{\textsf{fix}}}

\newcommand{\area}{{\textsf{area}}}

\newcommand{\W}{{\mathcal{W}}}

\newcommand{\M}{{\mathcal{M}}}

\title{Bijective Enumeration and Sign-Imbalance for\\ Permutation Depth and Excedances}
\author{Sen-Peng Eu
\institute{Department of Mathematics\\ National Taiwan Normal University\\Taiwan, ROC}
\institute{Chinese Air Force Academy\\
Taiwan, ROC}
\email{speu@math.ntnu.edu.tw}
\and
Tung-Shan Fu \qquad\qquad Yuan-Hsun Lo
\institute{Department of Applied Mathematics\\
National Pingtung University\\ Taiwan, ROC}
\email{\quad tsfu@mail.nptu.edu.tw \quad\qquad yhlo@mail.nptu.edu.tw}
}
\def\titlerunning{Sign-Imbalance for Permutation Depth and Excedances}
\def\authorrunning{S.-P. Eu, T.-S. Fu \& Y.-H. Lo}
\begin{document}
\maketitle

\begin{abstract}
We present a simplified variant of Biane's bijection between permutations and 3-colored Motzkin paths with weight that keeps track of the inversion number, excedance number and a statistic so-called depth of a permutation. This generalizes a result by Guay-Paquet and Petersen about a continued fraction of the generating function for depth on the symmetric group $\mathfrak{S}_n$ of permutations. In terms of weighted Motzkin path, we establish an involution on $\mathfrak{S}_n$ that reverses the parities of depth and  excedance numbers simultaneously, which proves that the numbers of permutations with even and odd depth (excedance numbers, respectively) are equal if $n$ is even and differ by the tangent number if $n$ is odd. Moreover, we present some interesting sign-imbalance results on permutations and derangements, refined with respect to depth and excedance numbers.
\end{abstract}

\section{Introduction}
\subsection{Preliminaries and background} 
Let $\mathfrak{S}_n$ be the symmetric group of permutations on $[n]:=\{1,2,\dots,n\}$. Given a permutation $\sigma=\sigma(1)\cdots \sigma(n)\in\mathfrak{S}_n$, an \emph{inversion} of $\sigma$ is a pair $(i,j)$ of indices such that $1\le i<j\le n$ and $\sigma(i)>\sigma(j)$. The \emph{inversion number} of $\sigma$, denoted by $\inv(\sigma)$, is defined to be the number of inversions of $\sigma$. 
Petersen and Tenner \cite{PT} presented a statistic called \emph{depth} for a Coxeter group, which is defined in terms of factorizations of the elements into reflections. In the case of symmetric group, the reflections of $\mathfrak{S}_n$ are the transpositions  $(i\,j)$, $1\le i<j\le n$, and the depth of a permutation $\sigma\in\mathfrak{S}_n$, denoted by $\depth(\sigma)$, is defined as
\begin{equation}
\depth(\sigma)=\min \left\{\sum_{r=1}^{k} (j_r-i_r) : \sigma=(i_1\, j_1)(i_2\, j_2)\cdots (i_k\, j_k) \right\}.
\end{equation}
They obtained a simple formula for the calculation of depth
\begin{equation}
\depth(\sigma)=\sum_{\sigma(i)>i} \big( \sigma(i)-i \big),
\end{equation}
which turns out to be one half of the \emph{total displacement} of $\sigma$  \cite[Problem~5.1.1.28]{Knuth} (also called \emph{Spearman's disarray} for $\mathfrak{S}_n$ in \cite{DG}).

For any two sequences $\{\gamma_h\}_{h\ge 0}$ and $\{\lambda_h\}_{h\ge 1}$, let $\mathfrak{F}(\gamma_h,\lambda_h)$ denote the series in $z$ defined by
\begin{equation} \label{eqn:J-fraction}
\mathfrak{F}(\gamma_h,\lambda_h)=\dfrac{1}{1-\gamma_0 z-\dfrac{\lambda_1 z^2}{1-\gamma_1 z-\dfrac{\lambda_2 z^2}{1-\gamma_2 z-\dfrac{\lambda_3 z^2}{\ddots}}}}
\end{equation}
These are Jacobi-type continued fractions (also known as $J$-\emph{fractions}). Flajolet \cite{Flajolet} gave an combinatorial interpretation of (\ref{eqn:J-fraction}) as the generating function for weighted Motzkin paths.

A \emph{Motzkin path} of length $n$ is a lattice path  from the origin to the point $(n,0)$ staying weakly above the $x$-axis, using the \emph{up step} $(1,1)$, \emph{down step} $(1,-1)$ and \emph{horizontal step} $(1,0)$. Let $U$, $D$ and $H$ denote an up step, a down step and a horizontal step, respectively. Let $\M_n$ denote the set of Motzkin paths of length $n$. 

For a Motzkin path $\mu=x_1\cdots x_n\in\M_n$, the \emph{height} of each step $x_j$ is the maximum $y$-coordinate achieved by $x_j$. For $x\in\{U,D,H\}$ and a nonnegative integer $h$, let $x^{(h)}$ denote a step $x$ at height $h$. The \emph{weight} $\omega(\mu)$ of $\mu$ is defined to be the product of the weight $\omega(x_j)$ of each step $x_j$ for all $j\in [n]$. 
By the \emph{area} under $\mu$, denoted by $\area(\mu)$, we mean the area of the region enclosed by $\mu$ and the $x$-axis.
Guay-Paquet and Petersen \cite{GP} defined a weight of $\mu$ by setting 
\begin{equation*}
\omega(x_i)=\begin{cases}
ht^{(2h-1)/2} &\mbox{if $x_i=U^{(h)}$ or $D^{(h)}$;} \\[0.4ex]
(2h+1)t^h    &\mbox{if $x_i=H^{(h)}$} 
\end{cases}
\end{equation*}
and
established a surjective map $\Phi:\sigma\mapsto\mu$ of $\mathfrak{S}_n$ onto $\M_n$ with
$\depth(\sigma)=\area(\mu)$, which yields the following  generating function for permutations with respect to depth:
\begin{equation} \label{eqn:expansion-depth-CF}
\sum_{n\ge 0}\sum_{\sigma\in\mathfrak{S}_n}  t^{\depth(\sigma)}z^n=\mathfrak{F}(\gamma_h,\lambda_h)
\end{equation}
with $\lambda_h=h^2 t^{2h-1}$ and $\gamma_h=(2h+1)t^h$, i.e.,
\begin{equation*} 
\sum_{n\ge 0}\sum_{\sigma\in\mathfrak{S}_n} t^{\depth(\sigma)}z^n=\dfrac{1}{1-z-\dfrac{tz^2}{1-3tz-\dfrac{4t^3z^2}{1-5t^2z-\dfrac{9t^5z^2}{1-7t^3z-\dfrac{16t^7z^2}{1-\cdots}}}}}.
\end{equation*} 
Guay-Paquet and Petersen \cite{GP} commented that map $\Phi:\mathfrak{S}_n\rightarrow\M_n$ is due to Foata and Zeilberger \cite{FZ} and that the statistic $\depth$ coincides with a statistic, $\mathsf{Edif}$, studied by Clarke, Steingr\'imsson and Zeng \cite{CSZ}.

Following \cite{SZ}, a \emph{3-colored Motzkin path} is a Motzkin path with three kinds of horizontal steps, denoted by $H_1$, $H_2$ and $H_3$. Using a double labeling scheme,
Biane \cite{Biane} established a bijection between permutations and 3-colored Motzkin paths with weight, which keeps track of the number of inversions. We refer to an alternative version of this bijection given by Sokal and Zeng in \cite[Section 6.2]{SZ}. Elizalde \cite{Elizalde} used this bijection to study the joint distribution of multiple statistics on $\mathfrak{S}_n$. 
Using Biane's method, we propose a multivariate weight function on the steps of a 3-colored Motzkin path and turn the surjective map $\Phi:\mathfrak{S}_n\rightarrow\M_n$ of Guay-Paquet and Petersen into a bijection that gives a generalization of (\ref{eqn:expansion-depth-CF}). (This resolves a conjecture raised by Petersen in an unpublished note \cite{Petersen}.) 

\subsection{Main results}
Let $\W_n$ denote the set of 3-colored Motzkin paths of length $n$ with a weight function $\omega$ defined by setting $\omega(H_3^{(0)})=p$ and 
\begin{equation} \label{eqn:step-weight}
\left\{
\begin{aligned}
\omega(U^{(h)}) &\in\{st^{2h-1}, st^{2h-1}q, \dots, st^{2h-1}q^{h-1}\}; \\
\omega(D^{(h)}) &\in\{q^{2h-1}, q^{2h}, \dots, q^{3h-2}\}; \\
\omega(H^{(h)}_1) &\in\{st^hq^h, st^hq^{h+1}, \dots, st^hq^{2h-1}\}; \\
\omega(H^{(h)}_2) &\in\{t^hq^h, t^hq^{h+1}, \dots, t^hq^{2h-1}\}; \\
\omega(H^{(h)}_3) &= pt^hq^{2h} 
\end{aligned}
\right.
\end{equation}  
for all $h\ge 1$.

Let $\sigma=\sigma(1)\cdots\sigma(n)\in\mathfrak{S}_n$. We say that an index $i\in [n]$ is an \emph{excedance} if $\sigma(i)>i$ and a \emph{fixed point} if $\sigma(i)=i$. Let $\exc(\sigma)$ and $\fix(\sigma)$ denote the number of excedances and fixed points of $\sigma$, respectively. For all positive integers $k$, we use the notations for $q$-integer $[k]_q=1+q+\cdots+q^{k-1}$ and  $[0]_q=0$. One of our main results is the following bijection.

\begin{thm} \label{thm:weight-preserving bijection}
There is a bijection $\Psi:\mathfrak{S}_n\rightarrow\W_n$ such that a permutation $\sigma\in\mathfrak{S}_n$ with $i$ inversions, $j$ fixed points, $k$ excedances and a depth of $\ell$ is carried to a 3-colored Motzkin path $\mu\in\W_n$ with weight $\omega(\mu)=q^{i} p^j s^{k} t^{\ell}$. Therefore, we have

\begin{equation} \label{eqn:refined-statistic-CF}
\sum_{n\ge 0}\left(\sum_{\sigma\in\mathfrak{S}_n} q^{\inv(\sigma)} p^{\fix(\sigma)} s^{\exc(\sigma)} t^{\depth(\sigma)}   \right)z^n=\mathfrak{F}(\gamma_h,\lambda_h)
\end{equation}  
with $\lambda_{h} = s [h]_q^2 (qt)^{2h-1}$ and $\gamma_h = \big((1+s)[h]_q+pq^h\big)(qt)^h$.
\end{thm}

\medskip
\noindent
{\bf Remarks.} 
By the weight given in (\ref{eqn:step-weight}), notice that a horizontal step $x\in\{H^{(h)}_2, H^{(h)}_3\}$ can be distinguished by the $q$-factor of $\omega(x)$ for each $h\ge 1$. That is, if $\omega(x)$ contains $q^{h+d}$ then $x=H^{(h)}_2$ ($H^{(h)}_3$, respectively) if $0\le d\le h-1$ ($d=h$, respectively). Thus, the bijection $\Psi:\mathfrak{S}_n\rightarrow\W_n$ can still be established with $H^{(h)}_2, H^{(h)}_3$ combined and the parameter $p$ omitted.

\medskip
The classical \emph{Euler numbers} $E_n$, defined by
\begin{align*}
\sum_{n\ge 0}E_n\frac{z^n}{n!} &=\tan z+\sec z \\
&= 1+z+\frac{z^2}{2!}+2\frac{z^3}{3!}+5\frac{z^4}{4!}+16\frac{z^5}{5!}+61\frac{z^6}{6!}+272\frac{z^7}{7!}+1385\frac{z^8}{8!}+\cdots,
\end{align*}
count the the number of \emph{alternating permutations} in $\mathfrak{S}_n$, i.e., $\sigma\in\mathfrak{S}_n$ such that $\sigma_1>\sigma_2<\sigma_3>\cdots \sigma_n$. The numbers $E_{2n}$ are called the \emph{secant numbers} and the numbers $E_{2n+1}$ are called the \emph{tangent numbers}.

Based on the bijection $\Psi:\mathfrak{S}_n\rightarrow\W_n$, we shall establish an involution on $\mathfrak{S}_n$ in terms of weighted Motzkin paths, which reverses the parities of depth and the number of excedances simultaneously. This proves that
the numbers of permutations with even and odd depth (excedance numbers, respectively) are equal if $n$ is even and differ by $E_n$ ($E_n$ up to sign, respectively) if $n$ is odd.

\begin{thm} \label{thm:depth-balance} There is an involution $\sigma\mapsto\sigma'$ on $\mathfrak{S}_n$ satisfying 
\begin{equation*}
\depth(\sigma)-\depth(\sigma')=\exc(\sigma)-\exc(\sigma')=\inv(\sigma)-\inv(\sigma')\in\{1,0,-1\}
\end{equation*}
and resulting in the following identities
\begin{enumerate}
\item ${\displaystyle
\sum_{\sigma\in\mathfrak{S}_{n}} (-1)^{\depth(\sigma)} = 
\begin{cases}
            E_n &  \mbox{for $n$ odd} \\
            0 & \mbox{for $n$ even.}
\end{cases}   
}$
\vspace{0.2cm}
\item ${\displaystyle
\sum_{\sigma\in\mathfrak{S}_{n}} (-1)^{\exc(\sigma)} = 
\begin{cases}
            (-1)^{\frac{n-1}{2}}E_n &  \mbox{for $n$ odd} \\
            0 & \mbox{for $n$ even.}  
\end{cases}
}$                   
\end{enumerate}

\end{thm}

\medskip
Note that the identity in Theorem \ref{thm:depth-balance}(ii) is a classical result of Euler.
When the fixed points of $\sigma$ are ignored, we present some interesting results on the sign imbalances of permutations and derangements, refined with respect to depth and excedance numbers. 

\begin{thm} \label{thm:sign-imbalance-permutations}
For $n\ge 1$, we have
\begin{equation}
\sum_{\sigma\in\mathfrak{S}_n} (-1)^{\inv(\sigma)} s^{\exc(\sigma)} t^{\depth(\sigma)}=(1-st)^{n-1}.
\end{equation}
\end{thm}

A \emph{derangement} of size $n$ is a permutation in $\mathfrak{S}_n$ that contains no fixed point. Let $\DDD_n\subset\mathfrak{S}_n$ be the set of derangements of size $n$. 
Define
\begin{equation}
F_n=F_n(s,t)=\sum_{\sigma\in\DDD_n} (-1)^{\inv(\sigma)} s^{\exc(\sigma)} t^{\depth(\sigma)}.
\end{equation}
Several of the initial polynomials $F_n(s,t)$ are listed below:
\begin{align*}
F_1(s,t) &= 0,\\
F_2(s,t) &= -st,\\
F_3(s,t) &= s(1+s)t^2, \\
F_4(s,t) &= s^2t^2-s(1+s)^2t^3, \\
F_5(s,t) &= -2s^2(1+s)t^3+s(1+s)^3t^4.
\end{align*}
Collected in powers of $t$, the coefficient $s$-polynomials of $F_n(s,t)$ for  $2\le n\le 9$ are listed in Table \ref{tab:step-by-step}. 
We obtain a neat expression for the sign imbalance of the joint distribution of depth and excedance numbers over $\DDD_n$.

\begin{table}[ht]
\caption{The $s$-polynomial coefficients of $F_n(s,t)$ in powers of $t$ for $2\le n\le 9$.}
\centering
{\small
\begin{tabular}{c|cccccccc}
      & $t$  &  $t^2$   &  $t^3$ & $t^4$ & $t^5$ & $t^6$ & $t^7$ & $t^8$  \\
\hline
$F_2$ & $-s$ &          &        &       &       &       &       &        \\
$F_3$ &      & $s(1+s)$ &        &       &       &       &       &        \\
$F_4$ &      &  $s^2$   &  $-s(1+s)^2$     &       &       &       &       &       \\
$F_5$ &      &          &  $-2s^2(1+s)$    &  $s(1+s)^3$   &       &       &       &       \\
$F_6$ &      &          &  $-s^3$          &  $3s^2(1+s)^2$ & $-s(1+s)^4$     &       &       &    \\
$F_7$ &      &          &                  &  $3s^3(1+s)$   & $-4s^2(1+s)^3$     & $s(1+s)^5$      &    &    \\
$F_8$ &      &          &                  &  $s^4$         & $-6s^3(1+s)^2$     & $5s^2(1+s)^4$   &  $-s(1+s)^6$     &       \\
$F_9$ &      &          &                  &                & $-4s^4(1+s)$     & $10s^3(1+s)^3$   &  $-6s^2(1+s)^5$     &  $s(1+s)^7$      
\end{tabular}
}
\label{tab:step-by-step}
\end{table}



\begin{thm} \label{thm:sign-imbalance-derangements} We have
\begin{equation*}
\sum_{n\ge 1} \left(\sum_{\sigma\in\DDD_n} (-1)^{\inv(\sigma)} s^{\exc(\sigma)} t^{\depth(\sigma)}\right)z^n=\sum_{k\ge 1} 
(-1)^k\left(\sum_{i=0}^{k-1}\binom{k-1}{i} s^{1+i}(1+s)^{k-1-i}z^{k+1+i} \right)t^k.
\end{equation*}
\end{thm}

\nocite{*}
\bibliographystyle{eptcs}
\bibliography{Bijective-Enumeration-References}
\end{document}